# ON OPTIMALITY OF THE BARRIER STRATEGY IN DE FINETTI'S DIVIDEND PROBLEM FOR SPECTRALLY NEGATIVE LÉVY PROCESSES

By R. L. Loeffen

*University of Bath*

We consider the classical optimal dividend control problem which was proposed by de Finetti [*Trans. XVth Internat. Congress Actuaries* **2** (1957) 433–443]. Recently Avram, Palmowski and Pistorius [*Ann. Appl. Probab.* **17** (2007) 156–180] studied the case when the risk process is modeled by a general spectrally negative Lévy process. We draw upon their results and give sufficient conditions under which the optimal strategy is of barrier type, thereby helping to explain the fact that this particular strategy is not optimal in general. As a consequence, we are able to extend considerably the class of processes for which the barrier strategy proves to be optimal.

**1. Introduction.** De Finetti [8] introduced the dividend model in risk theory. In this model the insurance company has the option to pay out dividends of its surplus to its beneficiaries up to the moment of ruin. De Finetti [8] argued that this should be done in an optimal way, namely such that the expected sum of the discounted paid out dividends from time zero until ruin is maximized. He proved that if the risk/surplus process evolves as a random walk with step sizes ±1, then an optimal way of paying out dividends is according to a barrier strategy, that is, there exists a constant $a^* \geq 0$, such that at each time epoch the excess of the net risk process over the level $a^*$ is paid out. In the case of continuous-time models, the problem of finding the optimal dividend strategy has been studied extensively in the Brownian motion setting [2, 21, 26, 31] and in the Cramér–Lundberg setting [4, 6, 14, 29], where by the former is meant that the risk process $X = \{X_t : t \geq 0\}$ is modeled by a Brownian motion plus drift and by the

latter that

$$X_t - X_0 = ct - \sum_{i=1}^{N_t} C_i,$$

where $C_1, C_2, \ldots$ are i.i.d. positive random variables representing the claims, $c > 0$ represents the premium rate and $N = \{N_t : t \geq 0\}$ is an independent Poisson process with arrival rate $\lambda$. Note that traditionally in the Cramér–Lundberg model it is assumed that $X$ drifts to infinity, but this condition is not necessary to formulate the problem. Very recently, Avram, Palmowski and Pistorius [3] considered the case where the risk process is given by a general spectrally negative Lévy process. Explanations for why this particular process serves as an appropriate generalization of the classical compound Poisson risk process can be found in, for example, [13, 18, 23]. It has been proved that in the Brownian motion setting and in the Cramér–Lundberg setting with exponentially distributed claims, an optimal dividend strategy is formed by a barrier strategy. No other explicit examples of spectrally negative Lévy processes have been given for which the same can be said. On the other hand, Azcue and Muler [4], Section 10.1, have found an example for which the optimal strategy is not a barrier strategy. Further, Avram, Palmowski and Pistorius [3] have given a sufficient condition involving the generator of the Lévy process for optimality of the barrier strategy. Besides finding the optimal strategy, a large body of literature exists [9, 12, 15, 16, 19, 23, 24, 25, 28, 32, 33] in which expressions are derived for, for example, the expected time of ruin, the moments of the expected paid out dividends and the Gerber–Shiu discounted penalty function, under the assumption that the insurance company pays out dividends according to a barrier strategy; the main motivation being the fact that the barrier strategy is optimal in (at least) the aforementioned two examples.

In this article motivated by the long history and broad interest of this control problem, we will shed new light on optimality of the barrier strategy when the risk process is modeled by a spectrally negative Lévy process. Using the setup and results from Avram, Palmowski and Pistorius [3], we show that the shape of the so-called scale functions of spectrally negative Lévy processes plays a central role. Further we will prove optimality of the barrier strategy if an easily checked analytical condition is imposed on the jump measure of the underlying Lévy process. This enables us to extend considerably the class of processes for which this strategy is optimal.

The outline of this paper is as follows. In Sections 2 and 3 we state the problem and briefly introduce scale functions. We present our main results in Section 4 and prove them in Section 5 using some earlier results from Avram, Palmowski and Pistorius [3]. We then conclude by giving some explicit examples to illustrate our results.

**2. Problem setting.** Let $X = \{X_t : t \geq 0\}$ be a spectrally negative Lévy process on a filtered probability space $(\Omega, \mathcal{F}, \mathbb{F} = \{\mathcal{F}_t : t \geq 0\}, \mathbb{P})$ satisfying the usual conditions. We denote by $\{\mathbb{P}_x, x \in \mathbb{R}\}$ the family of probability measures corresponding to a translation of $X$ such that $X_0 = x$, where we write $\mathbb{P} = \mathbb{P}_0$. Further $\mathbb{E}_x$ denotes the expectation with respect to $\mathbb{P}_x$ with $\mathbb{E}$ being used in the obvious way. Let the Lévy triplet of $X$ be given by $(\gamma, \sigma, \nu)$, where $\gamma \in \mathbb{R}$, $\sigma \geq 0$ and $\nu$ is a measure on $(0, \infty)$ satisfying

$$\int_{(0,\infty)} (1 \wedge x^2) \nu(dx) < \infty.$$

The Laplace exponent of $X$ is given by

$$\psi(\theta) = \log(\mathbb{E}(e^{\theta X_1})) = \gamma \theta + \tfrac{1}{2} \sigma^2 \theta^2 - \int_{(0,\infty)} (1 - e^{-\theta x} - \theta x \mathbf{1}_{\{0 < x < 1\}}) \nu(dx)$$

and is well defined for $\theta \geq 0$. Note that in the Cramér–Lundberg setting $\sigma = 0$, $\nu(dx) = \lambda F(dx)$ where $F$ is the law of $C_1$ and $\gamma = c - \int_{(0,1)} x \nu(dx)$. We exclude the case that $X$ has monotone paths. The process $X$ will represent the risk/surplus process of an insurance company before dividends are deducted.

We denote a dividend or control strategy by $\pi$, where $\pi = \{L_t^\pi : t \geq 0\}$ is a nondecreasing, left-continuous $\mathbb{F}$-adapted process which starts at zero. $L_t^\pi$ will represent the cumulative dividends the company has paid out until time $t$ under the control $\pi$. We define the controlled (net) risk process $U^\pi = \{U_t^\pi : t \geq 0\}$ by $U_t^\pi = X_t - L_t^\pi$. Let $\sigma^\pi = \inf\{t > 0 : U_t^\pi < 0\}$ be the ruin time and define the value function of a dividend strategy $\pi$ by

$$v_\pi(x) = \mathbb{E}_x \left[ \int_0^{\sigma^\pi} e^{-qt} \, dL_t^\pi \right],$$

where $q > 0$ is the discount rate. By definition it follows that $v_\pi(x) = 0$ for $x < 0$. A strategy $\pi$ is called admissible if ruin does not occur by a dividend payout, that is, $L_{t+}^\pi - L_t^\pi \leq U_t^\pi$ for $t < \sigma^\pi$. Let $\Pi$ be the set of all admissible dividend policies. The control problem consists of finding the optimal value function $v_*$ given by

$$v_*(x) = \sup_{\pi \in \Pi} v_\pi(x)$$

and an optimal strategy $\pi_* \in \Pi$ such that

$$v_{\pi_*}(x) = v_*(x) \qquad \text{for all } x \geq 0.$$

We denote by $\pi_a = \{L_t^a : t \geq 0\}$ the barrier strategy at level $a$ which is defined by $L_0^a = 0$ and

$$L_t^a = \left( \sup_{0 \leq s < t} X_s - a \right) \vee 0 \qquad \text{for } t > 0.$$

Note that $\pi_a \in \Pi$. Let $v_a$ denote the value function when using the dividend strategy $\pi_a$. In this paper we find sufficient conditions such that $v_*(x) = v_a(x)$ for all $x \geq 0$ for a certain specified $a$.

**3. Scale functions.** For each $q \geq 0$ there exists a function $W^{(q)} \colon \mathbb{R} \to [0, \infty)$, called the $(q$-)scale function of $X$, which satisfies $W^{(q)}(x) = 0$ for $x < 0$ and is characterized on $[0, \infty)$ as a strictly increasing and continuous function whose Laplace transform is given by

$$\int_0^\infty e^{-\theta x} W^{(q)}(x)\, dx = \frac{1}{\psi(\theta) - q} \qquad \text{for } \theta > \Phi(q),$$

where $\Phi(q) = \sup\{\theta \geq 0 \colon \psi(\theta) = q\}$ is the right-inverse of $\psi$. We write $W = W^{(0)}$. We will later on use the following relation:

$$(1) \qquad\qquad W^{(q)}(x) = e^{\Phi(q)x} W_{\Phi(q)}(x).$$

Here $W_{\Phi(q)}$ is the $(0$-)scale function of $X$ under the measure $\mathbb{P}^{\Phi(q)}$, where this measure is defined by the change of measure

$$\left.\frac{d\mathbb{P}^{\Phi(q)}}{d\mathbb{P}}\right|_{\mathcal{F}_t} = e^{\Phi(q)X_t - qt}.$$

The process $X$ under the measure $\mathbb{P}^{\Phi(q)}$ is still a spectrally negative Lévy process, but with a different Lévy triplet. In particular its Lévy measure is now given by $e^{-\Phi(q)x}\nu(dx)$. We refer to [22], Chapter 8, for more information on scale functions.

Throughout this paper we will use the term sufficiently smooth, whereby we mean the following. A function $f \colon \mathbb{R} \to \mathbb{R}$ which vanishes on $(-\infty, 0)$ is called sufficiently smooth at a point $x > 0$ if $f$ is continuously differentiable at $x$ when $X$ is of bounded variation and is twice continuously differentiable at $x$ when $X$ is of unbounded variation. A function is then called sufficiently smooth if it is sufficiently smooth at all $x > 0$; see [7] for conditions under which the scale function $W^{(q)}$ is sufficiently smooth. The derivative of $x \mapsto W^{(q)}(x)$ is denoted by $W^{(q)\prime}$.

Avram, Palmowski and Pistorius [3] showed that the value of the barrier strategy can be expressed in terms of scale functions in the following way.

PROPOSITION 1. *Assume $W^{(q)}$ is continuously differentiable on $(0, \infty)$. The value function of the barrier strategy at level $a \geq 0$ is given by*

$$v_a(x) = \begin{cases} \dfrac{W^{(q)}(x)}{W^{(q)\prime}(a)}, & \text{if } x \leq a, \\[2ex] x - a + \dfrac{W^{(q)}(a)}{W^{(q)\prime}(a)}, & \text{if } x > a. \end{cases}$$

The proof of Proposition 1 in [3] is based on excursion theory. An alternative proof where only basic fluctuation identities are used in conjunction with the strong Markov property is given in [28, 34]. Define now the (candidate) optimal barrier level by

$$a^* = \sup\{a \geq 0 : W^{(q)\prime}(a) \leq W^{(q)\prime}(x) \text{ for all } x \geq 0\},$$

where $W^{(q)\prime}(0)$ is understood to be equal to $\lim_{x \downarrow 0} W^{(q)\prime}(x)$. It follows that $a^* < \infty$ since $\lim_{x \to \infty} W^{(q)\prime}(x) = \infty$. Note that our definition of the optimal barrier level is slightly different than the one given by Avram, Palmowski and Pistorius [3]. It is easily seen that if an optimal strategy is formed by a barrier strategy, then the barrier strategy at $a^*$ has to be an optimal strategy.

**4. Main results.** We will now present the main results of this paper which give sufficient conditions for optimality of the barrier strategy $\pi_{a^*}$.

THEOREM 2. *Suppose $W^{(q)}$ is sufficiently smooth and*

$$(2) \qquad W^{(q)\prime}(a) \leq W^{(q)\prime}(b) \qquad \text{for all } a^* \leq a \leq b.$$

*Then the barrier strategy at $a^*$ is an optimal strategy.*

A drawback of condition (2) is that it involves the scale function for which closed form expressions are only known in a few cases. It would be better to have a condition which is directly given in terms of the Lévy triplet $(\gamma, \sigma, \nu)$ and the discount rate $q$. The second theorem entails exactly such a condition.

THEOREM 3. *Suppose that the Lévy measure $\nu$ of $X$ has a completely monotone density, that is, $\nu(dx) = \mu(x)\,dx$, where $\mu : (0, \infty) \to [0, \infty)$ has derivatives $\mu^{(n)}$ of all orders which satisfy*

$$(-1)^n \mu^{(n)}(x) \geq 0 \qquad \text{for } n = 0, 1, 2, \ldots.$$

*Then $W^{(q)\prime}$ is strictly convex on $(0, \infty)$ for all $q > 0$. Consequently, (2) holds and the barrier strategy at $a^*$ is an optimal strategy for the control problem.*

**5. Proof of main results.** Before proving the main results, we give two lemmas. Both lemmas are lifted from Avram, Palmowski and Pistorius [3]. We therefore do not give a proof of the first lemma which is a verification lemma involving a Hamilton–Jacobi–Bellman inequality. However, we do include a short proof of the second one as various arguments will be instructive to refer back to in the proof of Theorem 2.

Let $\Gamma$ be the operator acting on sufficiently smooth functions $f$, defined by

$$\Gamma f(x) = \gamma f'(x) + \frac{\sigma^2}{2} f''(x)$$
$$+ \int_{(0,\infty)} [f(x-y) - f(x) + f'(x)y\mathbf{1}_{\{0<y<1\}}]\nu(dy).$$

Lemma 4 (Verification lemma). *Suppose $\pi$ is an admissible dividend strategy such that $v_\pi$ is sufficiently smooth and for all $x > 0$*

$$\max\{\Gamma v_\pi(x) - q v_\pi(x), 1 - v_\pi'(x)\} \leq 0 \qquad \text{(HJB inequality).}$$

*Then $v_\pi(x) = v_*(x)$ for all $x \in \mathbb{R}$.*

Lemma 5. *Suppose $W^{(q)}$ is sufficiently smooth and suppose that*

$$(3) \qquad\qquad (\Gamma - q)v_{a^*}(x) \leq 0 \qquad \text{for } x > a^*.$$

*Then $v_{a^*}(x) = v_*(x)$ for all $x \in \mathbb{R}$.*

PROOF OF LEMMA 5. It suffices to show that under the conditions of Lemma 5 $v_{a^*}$ satisfies the conditions of the verification lemma. When $a^* = 0$ this is trivial because of (3), so we assume without loss of generality that $a^* > 0$. Because $W^{(q)}$ is sufficiently smooth and by Proposition 1, it follows that for any $a \geq 0$, $v_a(x)$ is sufficiently smooth at all $x \in (0,\infty)\backslash\{a\}$. By definition of $a^*$ and the assumed smoothness, we have $W^{(q)\prime\prime}(a^*) = 0$ when $X$ is of unbounded variation and hence $v_{a^*}(x)$ is also sufficiently smooth at $x = a^*$. Further $v_{a^*}'(x) \geq 1$ by definition of $a^*$. Since $(e^{-q(t\wedge\tau_0^-\wedge\tau_a^+)}W^{(q)}(X_{t\wedge\tau_0^-\wedge\tau_a^+}))_{t\geq 0}$ is a $\mathbb{P}_x$-martingale, one can deduce that

$$(4) \qquad\qquad (\Gamma - q)v_a(x) = 0 \qquad \text{for } 0 < x < a \text{ and } a > 0.$$

[Note that for $a \neq a^*$, $v_a(x)$ is not necessarily twice continuously differentiable in $x = a$ even if $W^{(q)\prime\prime}$ is continuous in $a$. Therefore $(\Gamma - q)v_a(x)$ is not necessarily continuous in $a$ and so (4) does not hold for $x = a$ in general.] In particular (4) holds for $a = a^*$. Hence together with (3), $v_{a^*}$ satisfies the HJB inequality. $\square$

PROOF OF THEOREM 2. First, we claim that

$$(5) \qquad\qquad \lim_{y\uparrow x}(\Gamma - q)(v_{a^*} - v_x)(y) \leq 0 \qquad \text{for } x > a^*.$$

We prove the claim for $X$ being of unbounded variation (the case of bounded variation is slightly easier). Let $x > a^*$. By assumption on the smoothness of the scale function, $v_x$ and $v_{a^*}$ are twice continuously differentiable on $(0, \infty)$, except for the possibility that $\lim_{y \uparrow x} v''_x(y) \neq \lim_{y \downarrow x} v''_x(y)$. We can use the dominated convergence theorem to deduce

$$\lim_{y \uparrow x} (\Gamma - q)(v_{a^*} - v_x)(y)$$

$$= \gamma (v'_{a^*} - v'_x)(x) + \frac{\sigma^2}{2} \left( v''_{a^*}(x) - \lim_{y \uparrow x} v''_x(y) \right) - q(v_{a^*} - v_x)(x)$$

$$+ \int_{(0, \infty)} \{ [(v_{a^*} - v_x)(x - z) - (v_{a^*} - v_x)(x)]$$

$$+ (v'_{a^*} - v'_x)(x) z \mathbf{1}_{\{0 < z < 1\}} \} \nu(dz).$$

Since we have by using Proposition 1:

(i) $\lim_{y \uparrow x} v''_x(y) \geq 0 = v''_{a^*}(x)$ where the inequality is by (2),

(ii) $(v'_{a^*} - v'_x)(u) \geq 0$ for $u \in [0, x]$, since for $u \in [0, a^*]$ $(v'_{a^*} - v'_x)(u) \geq 0$ by definition of $a^*$ and for $u \in (a^*, x]$ $(v'_{a^*} - v'_x)(u) \geq 0$ by (2); this implies that $(v_{a^*} - v_x)(x - z) \leq (v_{a^*} - v_x)(x)$ for all $z \geq 0$,

(iii) $(v_{a^*} - v_x)(x) \geq 0$ which follows from $v_{a^*}(a^*) \geq v_x(a^*)$ and (ii),

(iv) $v'_{a^*}(x) = v'_x(x) = 1$,

the claim follows.

We now prove by contradiction that (3) holds; the theorem is then proved by applying Lemma 5. Suppose there exist $x > a^* \geq 0$ such that $(\Gamma - q)v_{a^*}(x) > 0$. Then by (5) and the continuity of $(\Gamma - q)v_{a^*}$ we have $\lim_{y \uparrow x} (\Gamma - q)v_x(y) > 0$ which contradicts (4). $\square$

PROOF OF THEOREM 3. Since $\nu_{\Phi(q)}(dx) = e^{-\Phi(q)x} \mu(x) \, dx$ is the Lévy measure of the process $X$ under the measure $\mathbb{P}^{\Phi(q)}$, we have that $\nu_{\Phi(q)}(dx)$ has a completely monotone density, since the product of two completely monotone functions is completely monotone. It follows that $x \longmapsto \nu_{\Phi(q)}(x, \infty)$ is completely monotone, since $\frac{d}{dx} \nu_{\Phi(q)}(x, \infty) = -e^{-\Phi(q)x} \mu(x)$.

Let $\{\hat{H}_t : t \geq 0\}$ be the descending ladder height process of $X$. As $q > 0$, the process $X$ under $\mathbb{P}^{\Phi(q)}$ drifts to infinity and it follows that the process $\hat{H}$ under $\mathbb{P}^{\Phi(q)}$ (under a suitably chosen constant appearing in the local time at the minimum) is a killed subordinator with Lévy measure given by $\nu_{\Phi(q)}(x, \infty) \, dx$ (see, e.g., [22], Exercise 6.5). Hence the Lévy measure of $\hat{H}$ under $\mathbb{P}^{\Phi(q)}$ has a completely monotone density and consequently the Laplace exponent of $\hat{H}$ under $\mathbb{P}^{\Phi(q)}$ is a complete Bernstein function (see [20], Theorem 3.9.29). We may now use a result from Rao, Song and Vondraček [27], Theorem 2.3, combined with [30], Remark 2.2, to conclude that the renewal function of $\hat{H}$ under $\mathbb{P}^{\Phi(q)}$ defined by $\hat{U}_{\Phi(q)}(x) = \mathbb{E}^{\Phi(q)}(\int_0^\infty \mathbf{1}_{\{\hat{H}_t \in [0, x]\}} \, dt)$ has a completely monotone derivative.

It is well known that the scale function of a spectrally negative Lévy process which does not drift to minus infinity is equal (up to a multiplicative constant appearing in the local time) to the renewal function of the descending ladder height process (see, e.g., [5], Chapter VII.2). So we can say that $W_{\Phi(q)}(x) = \widehat{U}_{\Phi(q)}(x)$ and therefore $W'_{\Phi(q)}$ is completely monotone. A nonnegative function on $(0, \infty)$ with a completely monotone derivative is also known as a Bernstein function.

Because $W_{\Phi(q)}|_{(0,\infty)}$ is a Bernstein function, it admits the following representation, which is closely related to Bernstein's theorem, (see, e.g., [20], Chapter 3.9):

$$(6) \qquad W_{\Phi(q)}(x) = a + bx + \int_{(0,\infty)} (1 - e^{-xt}) \xi(dt), \qquad x > 0,$$

where $a, b \geq 0$ and $\xi$ is a measure on $(0, \infty)$ satisfying $\int_{(0,\infty)} (t \wedge 1) \xi(dt) < \infty$; in other words $W_{\Phi(q)}$ is the Laplace exponent of some (possibly killed) subordinator. From (6) and (1) it follows that

$$W^{(q)}(x) = e^{\Phi(q)x}(a + bx) + \int_{(0,\infty)} (e^{\Phi(q)x} - e^{-x(t-\Phi(q))}) \xi(dt).$$

By repeatedly using the dominated convergence theorem, we can now deduce

$$W^{(q)'''}(x) = f'''(x) + \int_{(0,\infty)} (\Phi(q)^3 e^{\Phi(q)x} + (t - \Phi(q))^3 e^{-x(t-\Phi(q))}) \xi(dt)$$

$$= f'''(x) + \int_{(0,\Phi(q)]} (\Phi(q)^3 e^{\Phi(q)x} - (\Phi(q) - t)^3 e^{(\Phi(q)-t)x}) \xi(dt)$$

$$+ \int_{(\Phi(q),\infty)} (\Phi(q)^3 e^{\Phi(q)x} + (t - \Phi(q))^3 e^{-x(t-\Phi(q))}) \xi(dt),$$

where $f(x) = e^{\Phi(q)x}(a + bx)$. Hence $W^{(q)'''}(x) > 0$ for all $x > 0$ and so $W^{(q)'}$ is strictly convex on $(0, \infty)$. Since $W^{(q)}$ is infinitely differentiable, we can now apply Theorem 2 to deduce that the barrier strategy at $a^*$ is optimal. □

## 6. Examples.

EXAMPLE FROM THEOREM 2. We now give an example to illustrate Theorem 2. Let $X$ be given by the Cramér–Lundberg model perturbed by Brownian motion, that is,

$$X_t = x + ct - \sum_{i=1}^{N_t} C_i + \sigma B_t,$$

where we let $C_1 \sim \text{Erlang}(2, \alpha)$ (i.e., sum of two independent exponentially random variables with parameter $\alpha$). Note that the Lévy measure $\nu(dx) =$ $\lambda \alpha^2 x e^{-\alpha x} dx$ (where $\lambda$ is the arrival rate of the Poisson process $\{N_t : t \geq 0\}$) does not have a completely monotone density. For this example a closed form expression for the $q$-scale function in terms of the roots of $\psi(u) = q$ can easily be found by inverting its Laplace transform by the method of partial fraction expansion. Indeed, we can write (for $q > 0$ and $\sigma > 0$)

$$\frac{1}{\psi(u) - q} = \frac{1}{cu - \lambda + (\lambda \alpha^2/(\alpha + u)^2) + (1/2)\sigma^2 u^2 - q} \times \frac{(\alpha + u)^2}{(\alpha + u)^2}$$

$$= \frac{(\alpha + u)^2}{(1/2)\sigma^2 \prod_{j=1}^4 (u - \theta_j)} = \sum_{j=1}^4 \frac{D_j}{u - \theta_j},$$

where $(\theta_j)_{j=1}^4$ are the (possibly complex) zeros (which are assumed to be distinct) of the polynomial $(\psi(u) - q)(\alpha + u)^2$ and $(D_j)_{j=1}^4$ are given by

$$D_j = \frac{1}{\psi(u) - q}(u - \theta_j)\Big|_{u=\theta_j} = \frac{(\alpha + \theta_j)^2}{(1/2)\sigma^2 \prod_{k=1, k \neq j}^4 (\theta_j - \theta_k)}.$$

The scale function is then given by

$$W^{(q)}(x) = \sum_{j=1}^4 D_j e^{\theta_j x} \qquad \text{for } x \geq 0.$$

We now choose the values of the parameters as follows: $c = 21.4$, $\lambda = 10$, $\alpha = 1$, $q = 0.1$ and for $\sigma$ we consider two cases, the case when $\sigma = 1.4$ and $\sigma = 2$. (For these choices of the parameter values, the zeros $(\theta_j)_{j=1}^4$ are indeed distinct.) Note that when $\sigma = 0$, this is exactly the example given by Azcue and Muler [4] for which the optimal strategy is not of barrier type. In the two figures the graphs of $W^{(q)\prime}$ and $(\Gamma - q)v_{a^*}(x)$ for the chosen parameters are plotted with the help of Matlab. When $\sigma = 1.4$, $a^* \approx 0.4$ and we see from Figure 1 that (2) and also (3) do not hold. When $\sigma = 2$ the minimum of the derivative has shifted; now $a^* \approx 10.5$ and we see from Figure 2 that (2) does hold. Consequently by (the proof of) Theorem 2, (3) must hold, which is confirmed by the figure.   $\square$

Examples from Theorem 3.   By Theorem 3, we have that when the Lévy measure is completely monotone, then the barrier strategy at $a^*$ is always an optimal strategy. There are many examples of spectrally negative Lévy processes which have such a feature and which have been used in the literature to model the risk process. We name as examples the $\alpha$-stable process which has Lévy density

$$\mu(x) = \lambda x^{-1-\alpha} \qquad \text{with } \lambda > 0 \text{ and } \alpha \in (0, 2)$$

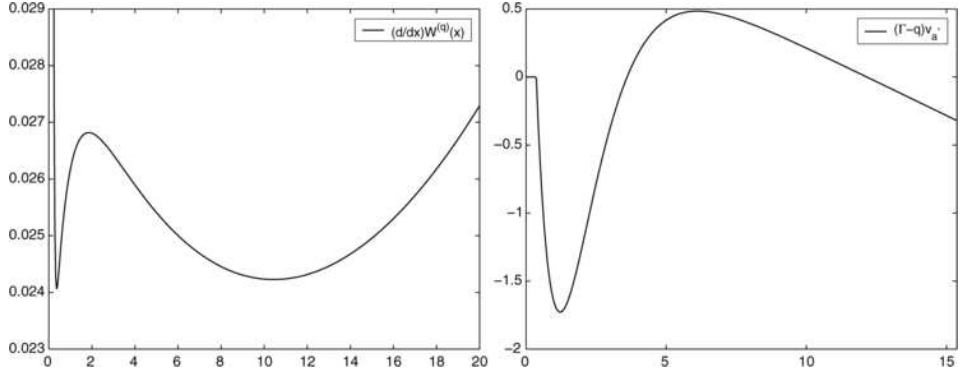

Fig. 1. $\sigma = 1.4$; left: $x \mapsto W^{(q)\prime}(x)$, right: $x \mapsto (\Gamma - q)v_{a^*}(x)$.

and is used in [13] and the (one-sided) tempered stable process which has Lévy density given by

$$\mu(x) = \lambda x^{-1-\alpha} e^{-\beta x} \qquad \text{with } \lambda, \beta > 0 \text{ and } -1 \le \alpha < 2.$$

The latter process includes other familiar Lévy processes, like the gamma process ($\alpha = 0$) which is considered in [11] and the inverse Gaussian process ($\alpha = 1/2$) which is used in [10] to model the risk process.

We can also conclude that the barrier strategy at $a^*$ is optimal, when we are in the Cramér–Lundberg setting where the claims have a distribution with a completely monotone probability density function. Some examples of these types of claim distributions which have been used in risk theory (see [1], Chapter I.2) are the heavy-tailed Weibull distribution

$$\mu(x) = crx^{r-1} e^{-cx^r} \qquad \text{with } c > 0 \text{ and } 0 < r < 1,$$

the Pareto distribution

$$\mu(x) = \alpha(1+x)^{-\alpha-1} \qquad \text{with } \alpha > 0$$

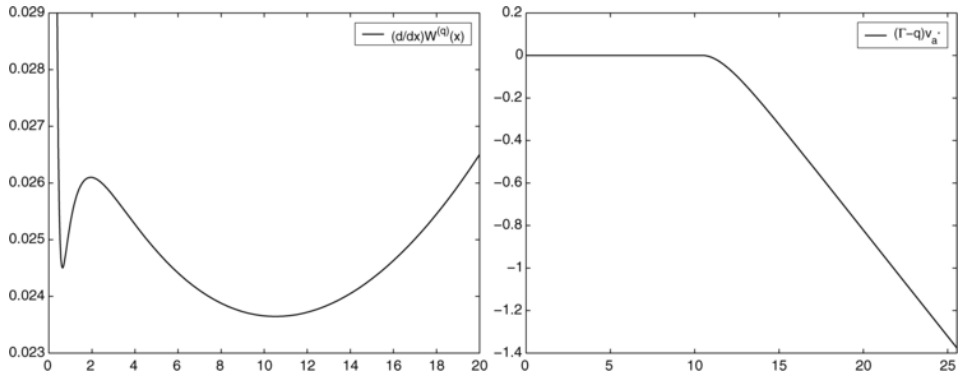

Fig. 2. $\sigma = 2$; left: $x \mapsto W^{(q)\prime}(x)$, right: $x \mapsto (\Gamma - q)v_{a^*}(x)$.

and the hyperexponential distribution

$$\mu(x) = \sum_{j=1}^{n} A_j \beta_j e^{-\beta_j x} \qquad \text{with } \beta_j, A_j > 0, j = 1, \ldots, n \text{ and } \sum_{j=1}^{n} A_j = 1.$$

Note that since in Theorem 3 there is no condition on the value of the Gaussian component $\sigma$, a barrier strategy will still form an optimal strategy if any one of the above examples is perturbed by Brownian motion.

For most spectrally negative Lévy processes an explicit expression for the $q$-scale function (and hence $a^*$) cannot be obtained. However, very recently Hubalek and Kyprianou [17] have found some new examples (including where the Lévy measure has a completely monotone density) for which the $q$-scale function is completely explicit.

## REFERENCES

[1] ASMUSSEN, S. (2000). *Ruin Probabilities. Advanced Series on Statistical Science & Applied Probability* **2**. World Scientific, River Edge, NJ. MR1794582

[2] ASMUSSEN, S. and TAKSAR, M. (1997). Controlled diffusion models for optimal dividend pay-out. *Insurance Math. Econom.* **20** 1–15. MR1466852

[3] AVRAM, F., PALMOWSKI, Z. and PISTORIUS, M. R. (2007). On the optimal dividend problem for a spectrally negative Lévy process. *Ann. Appl. Probab.* **17** 156–180. MR2292583

[4] AZCUE, P. and MULER, N. (2005). Optimal reinsurance and dividend distribution policies in the Cramér–Lundberg model. *Math. Finance* **15** 261–308. MR2132192

[5] BERTOIN, J. (1996). *Lévy Processes. Cambridge Tracts in Mathematics* **121**. Cambridge Univ. Press, Cambridge. MR1406564

[6] BÜHLMANN, H. (1970). *Mathematical Methods in Risk Theory. Die Grundlehren der mathematischen Wissenschaften* **172**. Springer, New York. MR0278448

[7] CHAN, T. and KYPRIANOU, A. E. (2009). Smoothness of scale functions for spectrally negative Lévy processes. Submitted.

[8] DE FINETTI, B. (1957). Su un'impostazion alternativa dell teoria collecttiva del rischio. In *Trans. XVth Internat. Congress Actuaries* **2** 433–443.

[9] DICKSON, D. C. M. and WATERS, H. R. (2004). Some optimal dividends problems. *Astin Bull.* **34** 49–74. MR2055537

[10] DUFRESNE, F. and GERBER, H. U. (1993). The probability of ruin for the inverse Gaussian and related processes. *Insurance Math. Econom.* **12** 9–22. MR1220358

[11] DUFRESNE, F., GERBER, H. U. and SHIU, E. S. W. (1991). Risk theory with the gamma process. *Astin Bull.* **21** 177–192.

[12] FROSTIG, E. (2005). The expected time to ruin in a risk process with constant barrier via martingales. *Insurance Math. Econom.* **37** 216–228. MR2172099

[13] FURRER, H. (1998). Risk processes perturbed by $\alpha$-stable Lévy motion. *Scand. Actuar. J.* **1** 59–74. MR1626676

[14] GERBER, H. U. (1969). Entscheidungskriterien für den zusammengesetzten Poisson-Prozess. *Mitt. Ver. Schweiz. Versich. Math.* **69** 185–227.

[15] GERBER, H. U., LIN, X. S. and YANG, H. (2006). A note on the dividends-penalty identity and the optimal dividend barrier. *Astin Bull.* **36** 489–503. MR2312676

[16] GERBER, H. U. and SHIU, E. S. W. (2004). Optimal dividends: Analysis with Brownian motion. *N. Am. Actuar. J.* **8** 1–20. MR2039869

[17] HUBALEK, F. and KYPRIANOU, A. E. (2007). Old and new examples of scale functions for spectrally negative Lévy processes. Preprint.

[18] HUZAK, M., PERMAN, M., ŠIKIĆ, H. and VONDRAČEK, Z. (2004). Ruin probabilities and decompositions for general perturbed risk processes. *Ann. Appl. Probab.* **14** 1378–1397. MR2071427

[19] IRBÄCK, J. (2003). Asymptotic theory for a risk process with a high dividend barrier. *Scand. Actuar. J.* **2** 97–118. MR1973214

[20] JACOB, N. (2001). *Pseudo Differential Operators and Markov Processes, Vol. I: Fourier Analysis and Semigroups.* Imperial College Press, London. MR1873235

[21] ZHANBLAN-PIKE, M. and SHIRYAEV, A. N. (1995). Optimization of the flow of dividends. *Uspekhi Mat. Nauk* **50** 25–46. MR1339263

[22] KYPRIANOU, A. E. (2006). *Introductory Lectures on Fluctuations of Lévy Processes with Applications.* Springer, Berlin. MR2250061

[23] KYPRIANOU, A. E. and PALMOWSKI, Z. (2007). Distributional study of de Finetti's dividend problem for a general Lévy insurance risk process. *J. Appl. Probab.* **44** 428–443. MR2340209

[24] LI, S. (2006). The distribution of the dividend payments in the compound Poisson risk model perturbed by diffusion. *Scand. Actuar. J.* **2** 73–85. MR2328678

[25] LIN, X. S., WILLMOT, G. E. and DREKIC, S. (2003). The classical risk model with a constant dividend barrier: Analysis of the Gerber–Shiu discounted penalty function. *Insurance Math. Econom.* **33** 551–566. MR2021233

[26] RADNER, R. and SHEPP, L. (1996). Risk vs. profit potential: A model for corporate strategy. *J. Econom. Dynamics Control* **20** 1373–1393.

[27] RAO, M., SONG, R. and VONDRAČEK, Z. (2006). Green function estimates and Harnack inequality for subordinate Brownian motions. *Potential Anal.* **25** 1–27. MR2238934

[28] RENAUD, J.-F. and ZHOU, X. (2007). Distribution of the present value of dividend payments in a Lévy risk model. *J. Appl. Probab.* **44** 420–427. MR2340208

[29] SCHMIDLI, H. (2006). Optimisation in non-life insurance. *Stoch. Models* **22** 689–722. MR2263862

[30] SONG, R. and VONDRAČEK, Z. (2006). Potential theory of special subordinators and subordinate killed stable processes. *J. Theoret. Probab.* **19** 817–847. MR2279605

[31] TAKSAR, M. I. (2000). Optimal risk and dividend distribution control models for an insurance company. *Math. Methods Oper. Res.* **51** 1–42. MR1742395

[32] YUEN, K. C., WANG, G. and LI, W. K. (2007). The Gerber–Shiu expected discounted penalty function for risk processes with interest and a constant dividend barrier. *Insurance Math. Econom.* **40** 104–112. MR2286657

[33] ZHOU, X. (2005). On a classical risk model with a constant dividend barrier. *N. Am. Actuar. J.* **9** 95–108. MR2211907

[34] ZHOU, X. (2006). Discussion of "On optimal dividend strategies in the compound Poisson model," by H. Gerber and E. Shiu. *N. Am. Actuar. J.* **10** 79–84.

DEPARTMENT OF MATHEMATICAL SCIENCES
UNIVERSITY OF BATH
CLAVERTON DOWN, BATH BA2 7AY
UNITED KINGDOM
E-MAIL: rll22@maths.bath.ac.uk